\pdfoutput=1
\documentclass{amsart}
\usepackage{a4wide}
\usepackage{amsmath}
\usepackage{amssymb}
\usepackage{graphicx,color}
\usepackage{hyperref}
\theoremstyle{definition}
\newtheorem{definition}{Definition}

\theoremstyle{theorem}
\newtheorem{theorem}{Theorem}

\newtheorem{corollary}{Corollary}
\newtheorem{fact}{Fact}

\theoremstyle{remark}
\newtheorem{remark}{Remark}

\newcommand{\ri}{\operatorname{RI}}

\begin{document}
\title[Automatic computation of crosscap number]{Automatic computation of crosscap number of alternating knots}
\author{Kaito Yamada}
\address{866 Nakane, Hitachinaka, Ibaraki, 312-8508, Japan}
\email{nanigasi.py@gmail.com}
\author{Noboru Ito}
\address{866 Nakane, Hitachinaka, Ibaraki, 312-8508, Japan}
\email{nito@gm.ibaraki-ct.ac.jp}
\date{March 17, 2023}
\maketitle
\begin{abstract}
We specify the computational complexity of crosscap numbers  of alternating knots by introducing an automatic computation. 
For an alternating knot $K$, let $\mathcal{E}$ be the number of edges of its diagram.  Then 
there exists a code such that the complexity of this computation of the crosscap number of $K$ is estimated by $O(\mathcal{E}^3)$.  
\end{abstract}
\section{Introduction}
The genus of a surface bounded by a knot $K$ is studied for about nine decades; remarkable progress has been made for the orientable case (the knot genus).  However, for non-orientable case (the crosscap number), there has not been a   \emph{general} approach to compute the crosscap number $C(K)$ of a knot $K$, whereas some pioneers, Hatcher-Thurston \cite{HatcherThurston1985}, Teragaito \cite{Teragaito2004}, Hirasawa-Teragaito \cite{HirasawaTeragaito2006}, and  Ichihara-Mizushima \cite{IchiharaMizushima2010}, give computations of crosscap numbers of certain families of knots.    

In \cite{ItoYamada2021}, the author KY introduces a new tool  as a Eulerian graph $G$ of a given alternating knot diagram $D$ to compute crosscap numbers automatically by computer aid.  The plan is as in the following list: 
\begin{enumerate}
\item If $G$ is obtained from  $D$, $G$ is uniquely determined.      \label{Y1}
\item Given $G$, the reduced graph is defined. \label{Y2}
\item For $n$, by inputting the initial data of the splice-unknotting number one, the data generates the list of alternating knot diagrams with splice-unknotting number $n$.   \label{Y3}
\end{enumerate}
The paper \cite{ItoYamada2021}  provides (\ref{Y1}) and (\ref{Y2}).  In this paper, the task (\ref{Y3}) is accomplished. 
\section{Crosscap number and splice-unknotting number}
\begin{definition}[non-orientable genus / crosscap number]
Let $K$ be a knot.  Then $C(K)$ is the minimum number among the first Betti numbers of non-orientable surfaces $\Sigma$ ($\subset \mathbb{R}^3$) with $\partial(\Sigma)$ $=$ $K$.  
\end{definition}
\begin{definition}[$S^-$]
Given a crossing on a link  diagram, there are the two possible ways to splice it.  One of them preserves orientations and the other does not; the symbol $S^-$ denotes the latter way.  
\end{definition}
\begin{definition}
Let $\ri^-$ be the first Reidemeister move resolving a single crossing.  
Let $D$ be a knot projection of a knot $K$.  The nonnegative number $u^-(D)$ is the minimum number of splices of type $S^-$ among any sequences of splices of types $S^-$ and $\ri^-$.   
\end{definition}
\begin{fact}[Ito-Takimura \cite{ItoTakimura2020}, Kindred \cite{ Kindred2020}]
Let $D$ be any prime alternating knot diagram of a knot $K$.  
\[
u^- (D) = C(K).  
\]  
\end{fact}
\begin{fact}[\cite{ItoTakimura2018}, cf.~ {\cite[twisted $S^+$ move/bridging operation]{ItoTakimura2022}}]\label{CorKD}
Every prime alternating knot diagram is obtained by applying bridging operation successively some times.   

More precisely, every prime alternating knot diagram $D$ with $u^-(D)=n$ is obtained from the knot projection with no crossing by applying bridging operation successively $n$ times.   
\end{fact}
\section{Main Result}   
\begin{theorem}\label{MainResult}
For an alternating knot $K$, let $C(K)$ be the crosscap number.  
A knot Eulerian graph is given by a knot diagram as in Definition~\ref{KEG}.  
Let $E$ be the number of edges of a knot Eulerian graph given by a knot diagram.  Then the complexity of the computation of $C(K)$ is bounded by $O(E^3)$.  
\end{theorem}
The number $V$ of vertices of a knot Eulerian graph is the less than or equal to the crossing number $\mathcal{V}$ of a knot diagram, i.e., $O(V) \le O(\mathcal{V})$.  Let $\mathcal{E}$ be the number of edges.  
$O(E)=O(V) \le O(\mathcal{V}) = O(\mathcal{E})$.   
\begin{corollary}
For an alternating knot $K$, let $C(K)$ be the crosscap number.  Let $\mathcal{E}$ be the number of edges of its diagram.    
Then there exists a code such that the complexity of the computation of $C(K)$ is bounded by $O(\mathcal{E}^3)$.
\end{corollary}
\begin{remark}
Even if $K$ is a non-alternating knot, the above estimation works in some case (cf.~\cite[Theorem~2, Fig.~1]{ItoTakimura2020}).  However,  the extent to which it will work is yet to be solved.  
\end{remark} 
\section{A short review of knot Eulerian graphs and update}\label{CAid}
\begin{definition}[twist region \cite{KalfagianniLee2016}]
A \emph{twist region} of a diagram consists of maximal collections of bigon regions arranged end to end; we suppose that the crossings in each twist region occur in an alternating fashion.    
\end{definition}
Note that a single crossing adjacent to no bigon is also a twist region.  
\begin{definition}[knot Eulerian graph]\label{KEG}
Every twisted region of an alternating knot diagram $D$ is given by a bridge operation increasing odd/even crossings, and $\ri^+$ if necessary.   The twisted region with odd/even crossings 
is represented as in Figure~\ref{Vertex}.  We call it an odd/even \emph{vertex} or a vertex simply.  
Every vertex has four endpoints that are two inputs and two outputs.  
Each vertex has an information as follows:
\begin{itemize}
\item By knot diagram, it is natural to suppose that every input corresponds to the unique output.  
\item  Valency of each vertex is exactly four.  
\end{itemize}
\begin{figure}[h!]
\includegraphics[width=12cm]{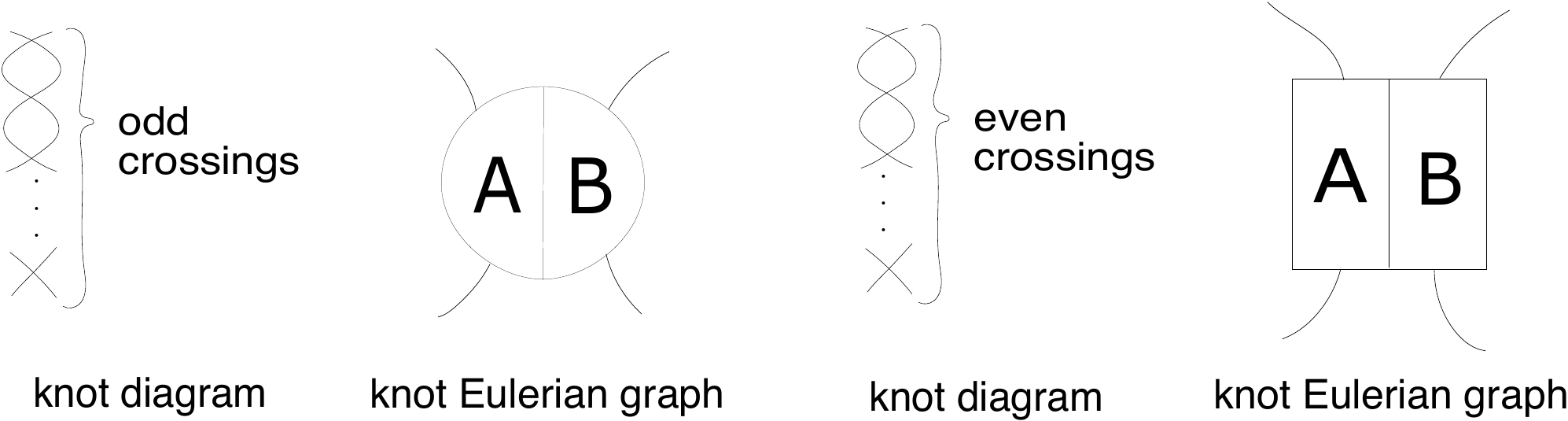}
\caption{Twisted regions of knot   diagrams and vertices of a knot Eulerian graph.}\label{Vertex}
\end{figure}
Further, we set the following rule for constructing our actual code.
\begin{itemize}
\item 
We often decompose each twisted region with even crossings into a single $\ri^+$ and a twisted region consisting of odd crossings.  
\item For every vertex, orientations of inputs or outputs are induced by knot diagrams as in Figure~\ref{Vertex}.  Each twist region with at least two crossings naturally gives a symmetric axis of the vertex and the left and right sides are called \emph{poles} and they are denoted by $A$ and $B$ in arbitrary way.   The  symmetric axis is called an \emph{axis}.  There are variations of axes by orientations (Figure~\ref{VertexVariation}).  
\begin{figure}[htbp]    
\centering
\includegraphics[width=10cm]{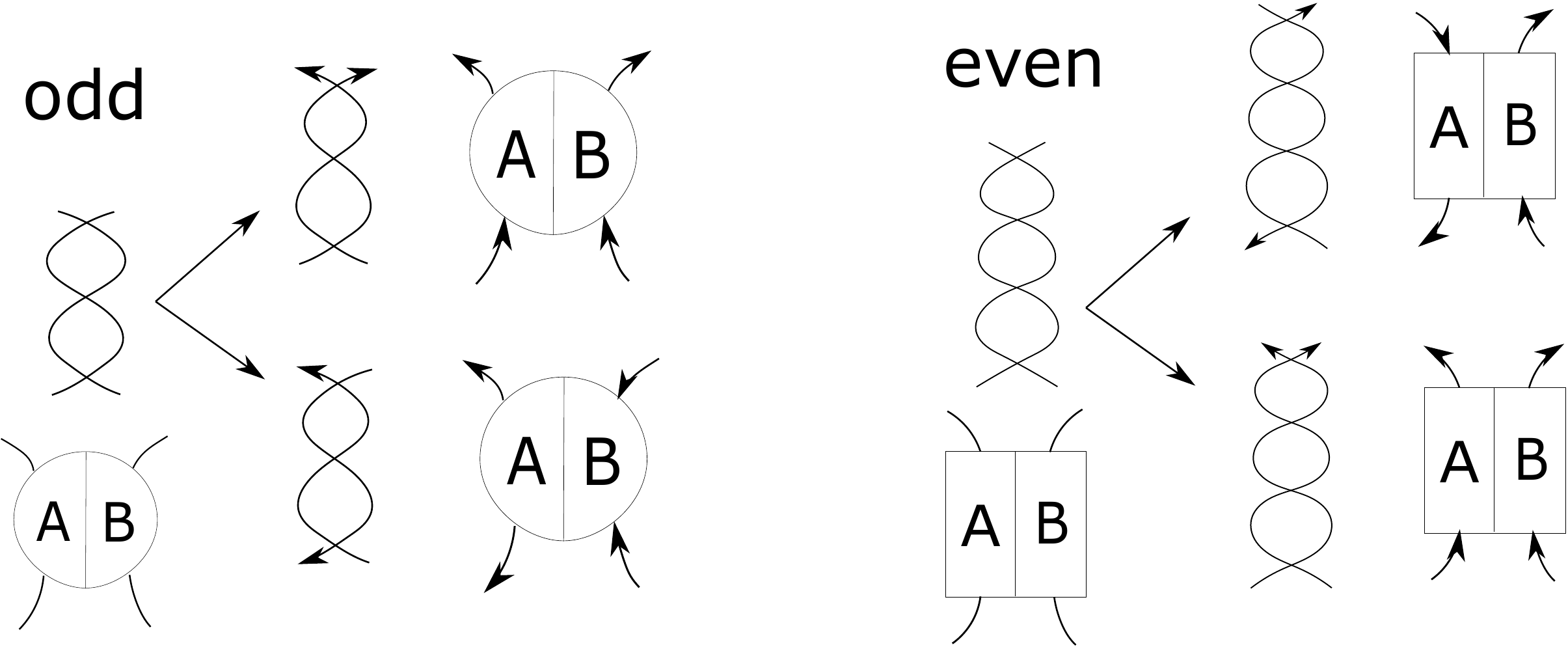} 
   \caption{Variations of vertices}
   \label{VertexVariation}
\end{figure}  
If a twist region has exactly a single crossing, we define the axis in arbitrary way.      
\end{itemize}
By the above definition, an alternating knot diagram implies vertexes; further the knot diagram has information of connections of the vertexes, which connections are canonically presented by edges. To code the edges, it is sufficient to use information of endpoints on poles.   Then a tuple $(u, v, P_u, P_v)$ denotes an edge, which satisfies conditions as follows.  
\begin{itemize}
\item (the first element of the tuple) $u$ denotes the vertex that is the starting point of the  edge.  
\item (the second element of the tuple) $v$ denotes the vertex that is the end point of the  edge.     
\item (the third element of the tuple) $P_u$ is the pole on which $u$ is.   
\item (the fourth element of the tuple) $P_v$ is the pole on which $v$ is.   
\end{itemize}

By the above construction, every alternating knot diagram gives a graph with four valences.   
The resulting graph is a Eulerian graph; we shall call it a \emph{knot Eulerian graph}.  
\end{definition}
\begin{remark}
In the previous paper \cite{ItoYamada2021}, $T_u$ denotes the vertex whereas we use $P_u$ that denotes it; $A$ and $B$ are called types whereas we call them \emph{poles}.   Note also that we do not need any empty vertex in this paper.  
\end{remark}
\begin{fact}[\cite{ItoYamada2021}]\label{PropKE}
An oriented knot projection gives a knot Eulerian graph.  
\end{fact}
In the rest of this paper, a vertex and an edge of a knot Eulerian graph are called a vertex and an edge, respectively.  However, to avoid any confusion, a vertex and an edge of a knot diagram is called a \emph{crossing} and a \emph{knot-edge}, respectively.  
\section{An estimation of computation complexity of bridge operations}\label{sec:Estimate}
Since it is possible to apply a bridge operation to \emph{every} knot-edge, we define two operations for knot Eulerian graph.  
\begin{itemize}
\item Decomposition of a vertex  (equivalently, a decomposition of a twist region of a knot diagram).  
\item $\ri^+$ (equivalently, the first Reidemeister move increasing a single crossing). 
\end{itemize} 
\subsection{Decomposition of a vertex}\label{sec:Deomp}
\begin{figure}[htbp] 
   \centering
   \includegraphics[width=5cm]{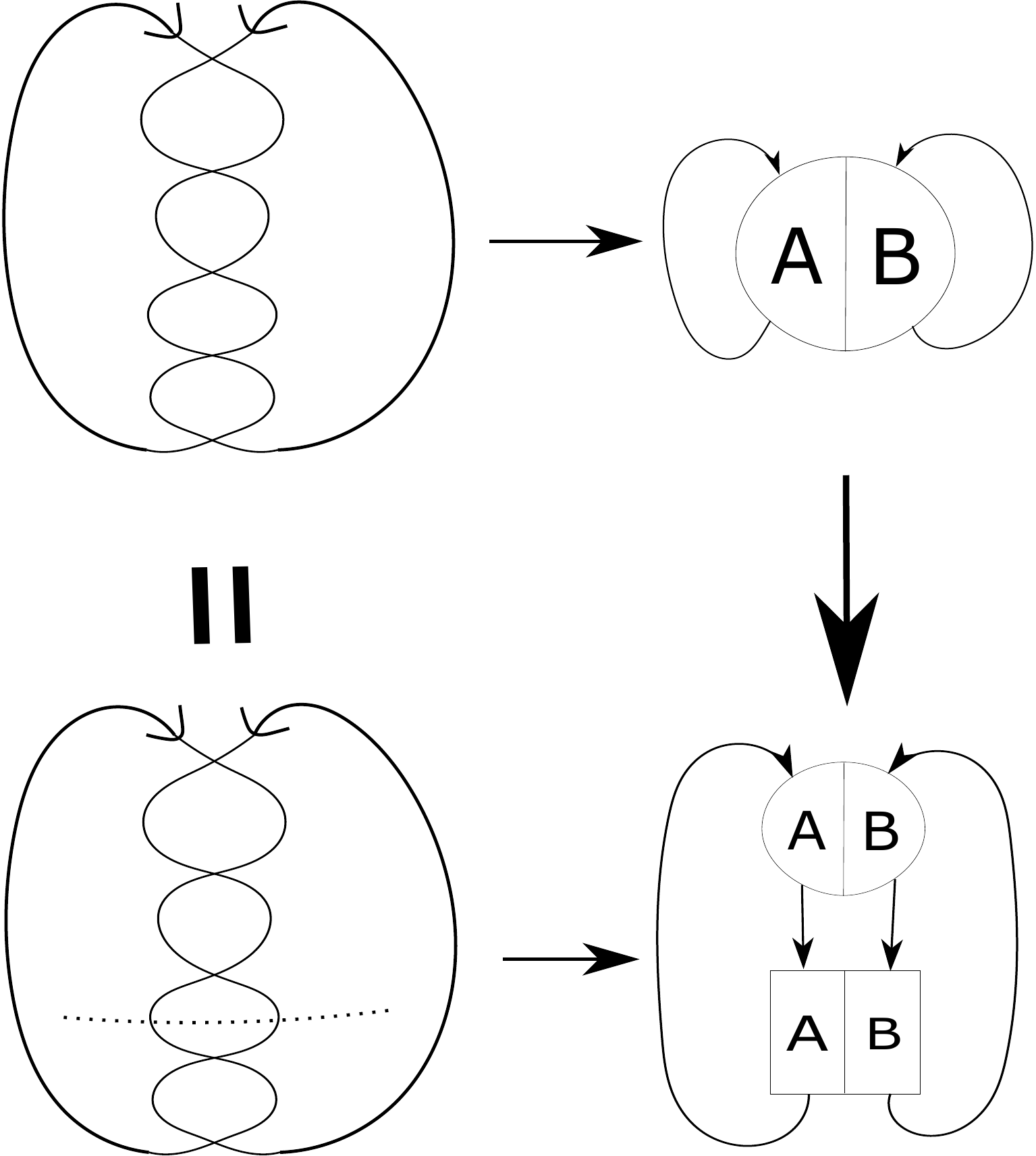} 
   \caption{An example of a  decomposition of a vertex}
   \label{Decomp}
\end{figure}

We define a decomposition of a vertex as in Figure~\ref{Decomp}.
There are two cases.  Let Odd (Even,~resp.) be the vertex type corresponding to a twist region of odd (even,~resp.) crossings.  
\begin{itemize}
\item Odd $\to$ Even and Odd / Odd and Even.  
\item Even $\to$ Odd and Odd  / Even and Even.  
\end{itemize}
\subsection{$\ri^+$}
\begin{itemize}
\item $(u, v, P_u, P_v)$ $\to$ $(u, O, P_u, A)$, $(O, O, B, A)$, and $(O, v, B, P_v)$
\end{itemize}
\subsection{Mixture of a decomposition of a vertex and applications of $\ri^+$}
\begin{figure}[htbp] 
   \centering
   \includegraphics[width=13cm]{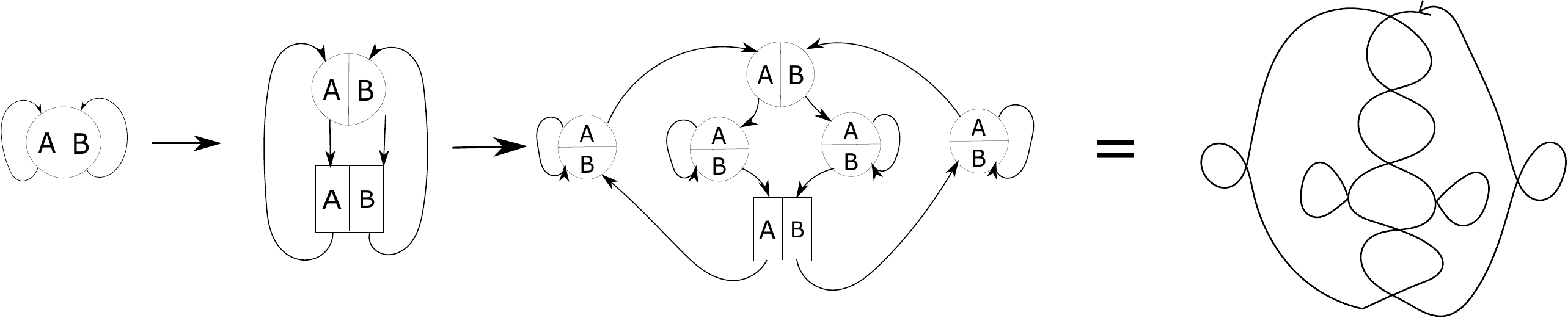} 
   \caption{An example of a mixture of a decomposition of a vertex and applications of $\ri^+$}
   \label{DeRI}
\end{figure}
Figure~\ref{DeRI} shows an example.  To make a computer program, we give a priority to a decomposition to fix the place where we apply a bridge operation, i.e., we always decompose a vertex before a $\ri^+$ is applied.     
\subsection{Selecting two edges which will be applied by a bridge operation}\label{Ecycle}
We estimate the complexity of computation of choice of the two places as follows: 
\begin{itemize}
\item We have a list of ordered edges $e_1, e_2, \dots, e_E$ that consists of a Euler cycle representing the knot Eulerian graph.   $\to$ $O(E)$.   
\item 
Two edges are chosen $\to$  $O(E^2)$.  
\end{itemize}
In the above process, these listing and selecting are totally estimated by $O(E^2)$.  
Let $e_s$ $=$ $(a, b, P_a, P_b)$ and $e_g$ $=$ $(c, d, P_c, P_d)$ ($s<g$) for two edges.
\subsection{Application of a bridge operation}\label{SecBridge}
\subsubsection{Algorithm}\label{SecAlg}
\begin{enumerate}
\item[Step~0:] We apply a decomposition as in Section~\ref{sec:Deomp} to vertices which connect to the selected above two edges.   
\item[Step~1:] We apply a single $\ri^+$ to  every edge.    
\item[Step~2:] Add a vertex $P$ with inputs $a, b$ and outputs $c, d$ as in Figure~\ref{Vertex}.  Define the four edges of the vertex $P$, i.e., $\alpha=(a, P, P_a, A)$, $\beta=(b, P, P_b, B)$, $\gamma=(P, c, B, P_c)$, $\delta=(P, d, A, P_d)$.  
\item[Step~3:] Add new edges $\{e'_{s+1}, e'_{s+2}, \dots, e'_{g-1}\}$ corresponding to $\{e_{s+1}, e_{s+2}, \dots, e_{g-1}\}$ where we give $e'_i$ by reversing the orientation of $e_i$.   
\item[Step~4:] Remove edges in $\{e_{s}, e_{s+1}, \dots, e_{g}\}$.     
\end{enumerate}
\begin{figure}[htbp]    \centering
\includegraphics[width=5cm]{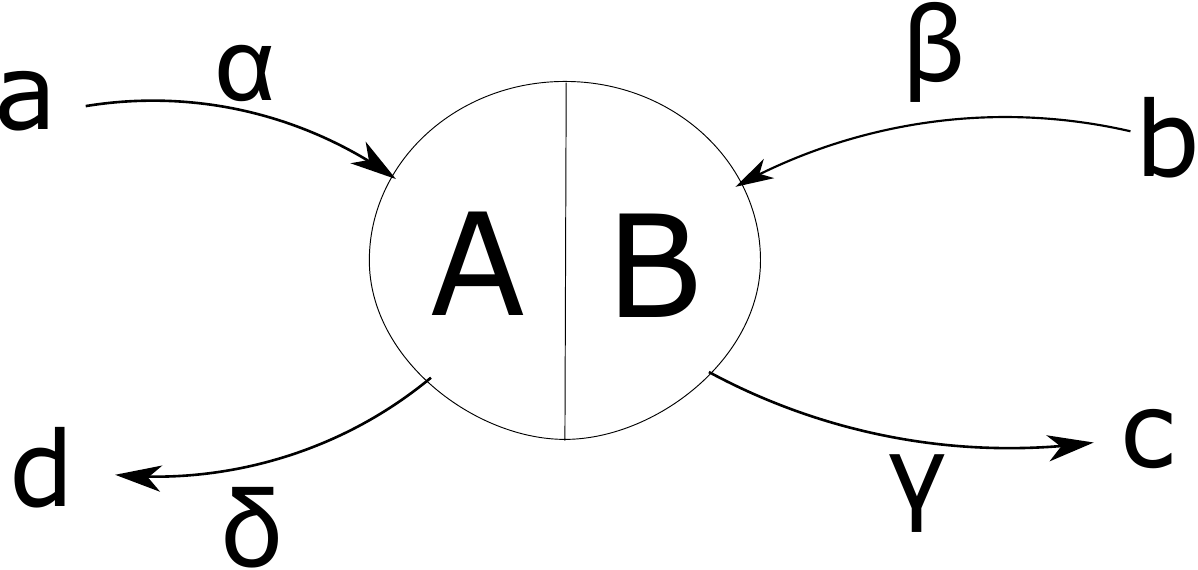} 
   \caption{Vertex $P$ in a general case}
   \label{AddingPUnori}
\end{figure}
\subsection{Examples: a bridge operation on a crosscap number one knot, Cases I--I\!I\!I}
\begin{figure}[h!]    
\centering
   \includegraphics[width=12cm]{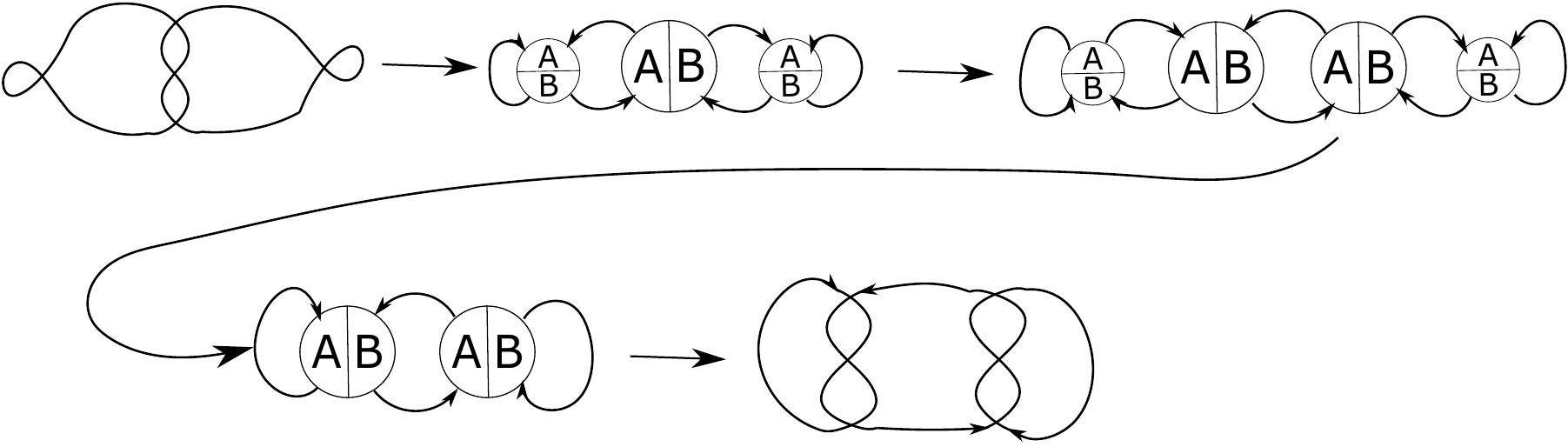} 
   \caption{Case~I}
   \label{CaseI}
\end{figure}

\begin{figure}[htbp]    
\centering
   \includegraphics[width=12cm]{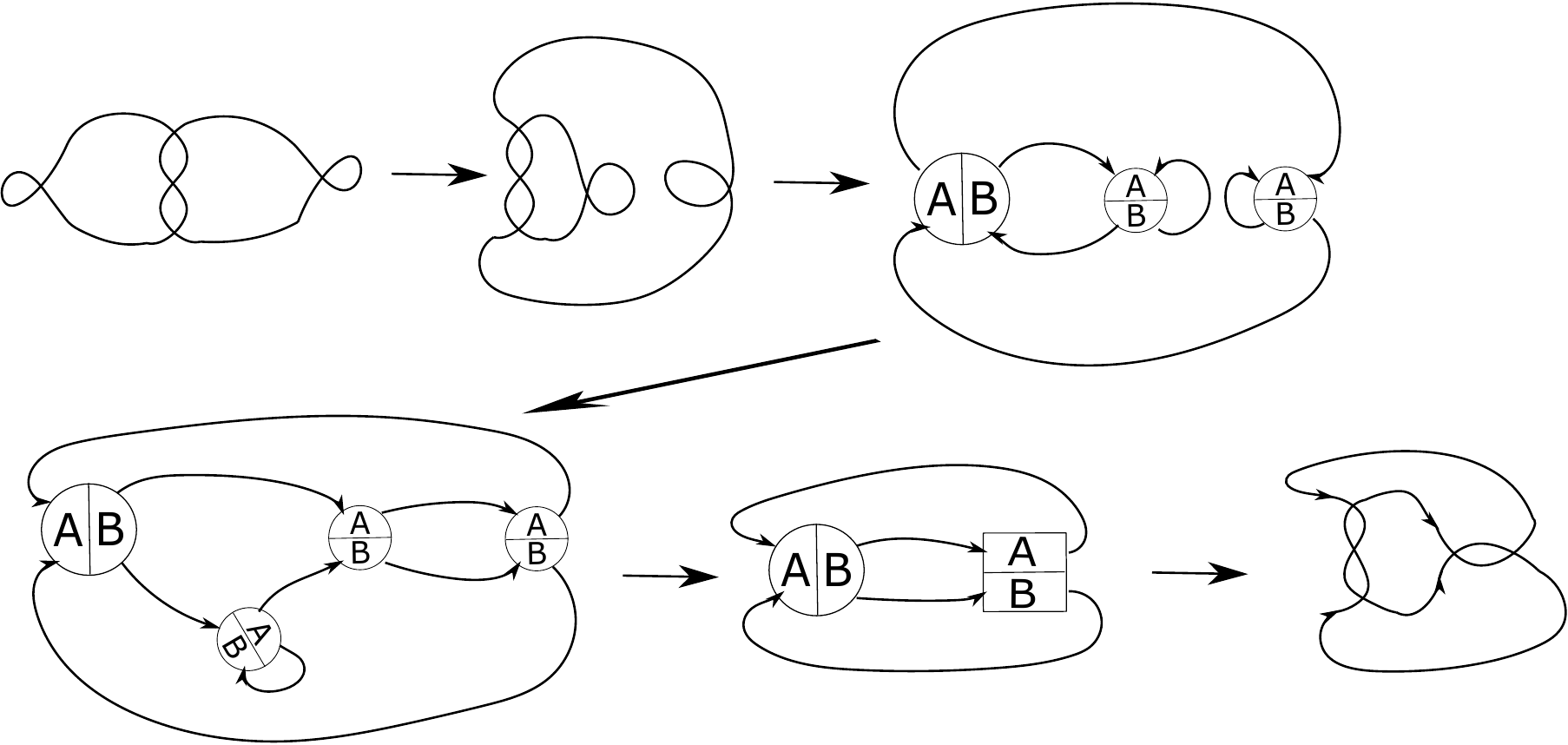} 
   \caption{Case~I\!I}
   \label{CaseII}
\end{figure}

\begin{figure}[htbp]    
\centering
   \includegraphics[width=12cm]{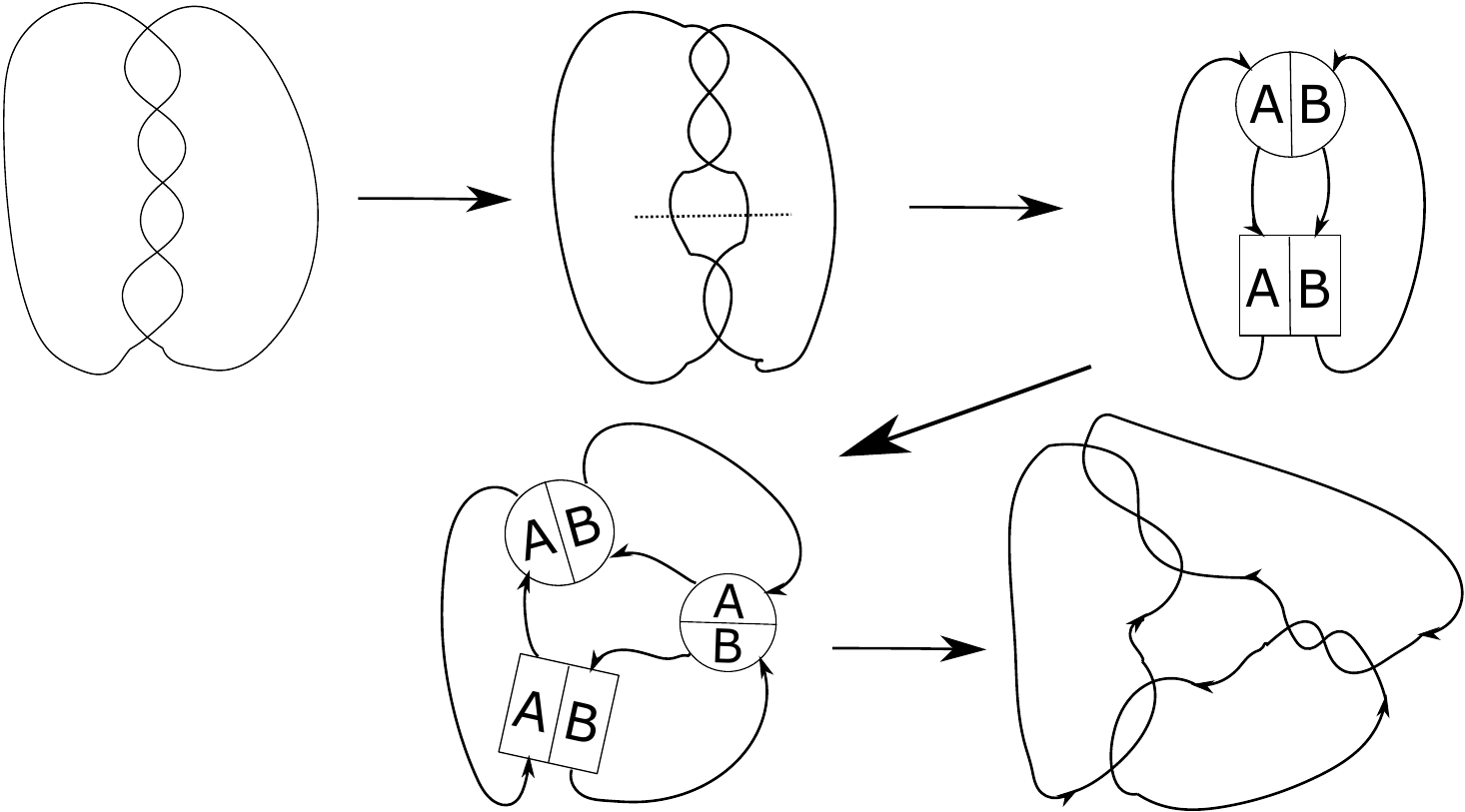} 
   \caption{Case~I\!I\!I}
   \label{CaseIII}
\end{figure}

We give examples: Cases I--I\!I\!I as in Figures~\ref{CaseI}--\ref{CaseIII}.  In Figures~\ref{CaseI}--\ref{CaseIII}, note that applying $\ri^+$'s at Step~1 often be omitted if they are wasted operations; in our actual code, we do not omit them.  
\subsection{An estimation of computation complexity of bridge operations}\label{Ebridge}
After applying steps in Section~\ref{Ecycle}, we give an estimation of the complexity of the computation in Section~\ref{SecBridge}.  
\begin{itemize}
\item We apply decompositions of four vertices that connect to the two selected edges $\to O(1)$.  
\item We apply $\ri^+$ for every edge $\to$ $O(E)$.  
\item We select two edges from decomposed edges including old/newborn edges by the above two steps (decompositions and $\ri^+$'s) $\to$ $O(1)$.  
\item We add a new vertex $P$  and new edges $\alpha$, $\beta$, $\gamma$, and $\delta$ $\to$ $O(1)$.   
\item We add new edges $e'_{s+1}$, \dots $e'_{g-1}$  $\to$ $O(E)$.   
\item We remove the edges from $e_s$ to $e_g$ $\to$ $O(E)$.  
\end{itemize}
The above process is totally estimated by $O(E)$.
\section{Proof of Theorem~\ref{MainResult}}
\begin{proof}
By Sections~\ref{Ecycle} and \ref{Ebridge}, for each fixing a pair which will be applied by a bridge operation, we complete a bridge operation in a code.  
\begin{itemize}
\item To fix each pair which will be applied by a bridge operation, the complexity of the computation is bounded by $O(E^2)$.  
\item To apply a bridge operation in a code, the complexity of the computation is bounded by $O(E)$.  
\end{itemize}
Hence $O(E^2) \times O(E)$ $=$ $O(E^3)$.  
\end{proof}
\section{Actual code}
For the detail, please see the actual code: 

\texttt{https://github.com/nanigasi-san/splus/blob/main/algorithms/main.py}

\section{Acknowledgements}
The work is partially supported by JSPS KAKENHI Grant Numbers  20K03604 and 22K03603 and Toyohashi Tech Project of Collaboration with KOSEN Grant Number 2309.  
\bibliographystyle{plain}
\bibliography{Ref}
\end{document}